# Pricing the Ramping Reserve and Capacity Reserve in Real Time Markets

Hongxing Ye, *Student Member, IEEE,* and Zuyi Li, *Senior Member, IEEE*

*Abstract*—The increasing penetration of renewable energy in recent years has led to more uncertainties in power systems. In order to maintain system reliability and security, electricity market operators need to keep certain reserves in the Security-Constrained Economic Dispatch (SCED) problems. A new concept, deliverable generation ramping reserve, is proposed in this paper. The prices of generation ramping reserves and generation capacity reserves are derived in the Affine Adjustable Robust Optimization framework. With the help of these prices, the valuable reserves can be identified among the available reserves. These prices provide crucial information on the values of reserve resources, which are critical for the long-term flexibility investment. The market equilibrium based on these prices is analyzed. Simulations on a 3-bus system and the IEEE 118-bus system are performed to illustrate the concept of ramping reserve price and capacity reserve price. The impacts of the reserve credit on market participants are discussed.

## Nomenclature

**Indices**

$i,k,l,t$    index for unit, uncertainty, line, and time
$m,n$    index for bus

**Functions and sets**

$\mathbb{E}(\cdot)$    expectation function
$\text{Tr}(\cdot), (\cdot)^\top, (\cdot)_j.$    trace, transpose, $j^{th}$ row of a matrix
$\Lambda(\cdot)$    cost as a function of adjustment matrix
$\Phi(\cdot)$    Lagrangian function
$(\cdot)^*$    optimal value of a variable
$\mathcal{G}(m)$    set of units located at bus m

**Constants**

$N_D, N_G, N_K$    number of buses, units and uncertainty constraints
$N_L, N_T$    number of lines and time intervals
$q_i, b_i$    cost coefficients for unit $i$
$\boldsymbol{Q_i}$    $\boldsymbol{Q_i} = q_i \boldsymbol{I_{N_T}} \in \mathbb{R}^{N_T \times N_T}$, $\boldsymbol{I_{N_T}}$ is unity matrix
$\boldsymbol{B_i}$    $\boldsymbol{B_i} = [b_i \cdots b_i]^\top \in \mathbb{R}^{N_T}$
$f_i$    cost related to unit commitment decision
$d_{m,t}, \boldsymbol{d}$    aggregated equivalent load demand on bus $m$ at $t$, $\boldsymbol{d} = [d_{1,t} \cdots d_{N_D,1} \cdots d_{N_D,N_T}]^\top$
$\bar{F}_l, \boldsymbol{F}$    branch flow limit, abstract vector $\boldsymbol{F} \in \mathbb{R}^{2N_L N_T}$
$\Gamma_{l,m}$    shift factor for line $l$ and bus $m$
$P_i^{\min}, P_i^{\max}$    minimum and maximum generation outputs
$\hat{I}_{i,t}, \hat{y}_{i,t}, \hat{z}_{i,t}$    unit on/off, start-up and shutdown indicators

The authors are with the Robert W. Galvin Center for Electricity Innovation at Illinois Institute of Technology, Chicago, IL 60616, USA. (e-mail: hye9@hawk.iit.edu; lizu@iit.edu).

$R_i^U, R_i^D$    unit ramping up/down limits
$\boldsymbol{D}$    time-load incidence matrix, $\boldsymbol{D} \in \mathbb{R}^{N_T \times N_D N_T}$
$\boldsymbol{A}, \boldsymbol{R_i}$    abstract matrix and vector for (1-3), $\boldsymbol{A} \in \mathbb{R}^{4N_T \times N_T}$, $\boldsymbol{R_i} \in \mathbb{R}^{4N_T}$.
$\boldsymbol{\Gamma_i}, \boldsymbol{\Gamma_d}$    abstract shift factor matrix for unit $i$ and load
$\boldsymbol{S}, \boldsymbol{h}$    polyhedron uncertainty set matrix and vector

**Variables**

$\boldsymbol{G_i}$    generation adjustment matrix, $\boldsymbol{G_i} \in \mathbb{R}^{N_T \times N_D N_T}$
$\boldsymbol{\pi}$    Prices. $\boldsymbol{\pi_i^e}$ is the energy LMP for unit $i$, $\boldsymbol{\pi_i^e} \in \mathbb{R}^{N_T}$; $\boldsymbol{\pi_i^r}$ is the reserve price vector for unit $i$, $\boldsymbol{\pi_i^r} \in \mathbb{R}^{4N_T}$; $\boldsymbol{\pi_i^a}$ is the marginal price for affine adjustment, $\boldsymbol{\pi_i^a} \in \mathbb{R}^{N_T \times N_D N_T}$.
$P_{i,t}$    generation output, $P_{i,t} \in \mathbb{R}$
$\boldsymbol{P_i}$    generation vector, $\boldsymbol{P_i} = [P_{i,1} \cdots P_{i,N_T}]^\top \in \mathbb{R}^{N_T}$
$\hat{\boldsymbol{P}}_i$    generation re-dispatch vector, $\hat{\boldsymbol{P}}_i \in \mathbb{R}^{N_T}$
$Q_{i,t}^{ru}, Q_{i,t}^{rd}$    generation ramping reserve, $Q_{i,t}^{ru} \in \mathbb{R}$, $Q_{i,t}^{rd} \in \mathbb{R}$
$Q_{i,t}^{cu}, Q_{i,t}^{cd}$    generation capacity reserve, $Q_{i,t}^{cu} \in \mathbb{R}$, $Q_{i,t}^{cd} \in \mathbb{R}$
$\boldsymbol{\epsilon}$    uncertainty vector, $\boldsymbol{\epsilon} \in \mathbb{R}^{N_T \times N_D}$
$\boldsymbol{\lambda}, \boldsymbol{\alpha_i}, \boldsymbol{\eta}$    Lagrangian multipliers for constraints (15, 17, 19)
$\boldsymbol{\gamma}, \boldsymbol{\beta_i}, \boldsymbol{\tau}$    Lagrangian multipliers for constraints (16,18,20)
$\Theta_i$    credits to uncertainty mitigators $i$

## I. Introduction

**T**HE renewable energy sources (RES), such as wind power generation, and price-sensitive demand response (DR) have attracted a lot of attentions recently. The total installed capacity of wind power in the U.S. reached 47 GW at the end of 2011 [1]. Several ISOs/RTOs, such as PJM, ISO New England, NYISO, and CAISO have initiated DR programs in their markets [2]. The essential objective to use renewable energy and initiate DR programs in electricity markets is to maximize the total social warfare as well as to protect the environment. However, they also pose new challenges to the system operators in electricity markets. Due to its intrinsic characteristics, the amount of available renewable energy is sometimes hard to predict. For instance, large-scale wind production varies from -20% to 20% of the installed wind capacity in Denmark on an hourly basis [3]. Prediction error for wind farms aggregated output by existing state-of-art method may fall into a range of 5%-20% of the total installed capacity [4]. In the meantime, the amount of un-predictable loads also increases in the wholesale market as the forecasting of price-sensitive loads relies on the forecasted price as input, which itself has significant error.

Surviving uncertainties is fundamentally important for the reliable and secure operation of a power system. If the system with pre-planned schedule cannot accommodate the deviation of the wind power and load from their forecasted values, the system operator may have to curtail wind energy or shed load in the real-time market (RTM). In order to keep a certain level of reliability or security, the ISO/RTO has to increase the ramping capability of the system to compensate the variations of wind energy and load in a short time [5], [6]. More efficient and reliable methods are required to determine the optimal reserves when the uncertainty level is high. Recently stochastic and robust approaches have been successfully applied by researchers to address the issues related to uncertainties in electricity markets. At the same time, market designers are also seeking effective market mechanism to address the uncertainty issues. For instance, intraday market (IDM) is now established between day-ahead market (DAM) and RTM in European countries since uncertainties on the intraday level are significantly smaller compared to those on the day-ahead level [7]. In the U.S., hour-ahead scheduling process (HASP) is employed by the California ISO [8].

Typical approaches to solving stochastic SCUC are scenario based [9], [10]. The basic idea is to generate enough samples for uncertain parameters with an assumption that their probability distribution function (PDF) is known. Those samples are then modeled in a Mixed Integer Linear/Quadratic Programming (MILP/MIQP) problem. The two main drawbacks of the scenario-based approaches are that the PDF is hard to obtain in some circumstances, and uncertainty accommodation is not guaranteed. In fact, the MILP/MIQP problem becomes intractable when the sample size is large. Comparing with stochastic optimization, two largest merits of robust optimization are that the solution can be immunized against all uncertainties and PDF is not required. In [11], [12], robust SCUC problem is solved in two stages. The first stage is to determine the unit commitment (UC) solution which is immunized against the worst case with the lowest cost. In the second stage, a feasible solution to SCED is obtained. Affinely adjustable robust optimization (AARO) models are proposed recently [13]–[15]. They employ an affine function to adjust the generation output following the load deviation. Recently, we propose an economic-efficient robust SCUC model with a fast solution approach [16].

Although applying robust techniques in SCUC/SCED receives a lot of attentions from researchers, it still remains a big challenge on how to credit the flexibilities in the U.S. electricity markets. In the existing Ancillary Service (AS) Market, the reserves are determined in advance [17]. The amount of required reserves is generally extracted from larger number of Monte Carlo simulations for the contingencies [18]. With the help of the AARO, those reserves are determined in one shot based on uncertainty information. Then a critical issue is how to price those reserves when there are no explicit reserve requirement constraints. On the other hand, some reserves are free byproducts in the co-optimization approach [17]. They are kept because the market participants want to get the energy profits. Moreover, some reserves are scarce resources due to their deliverability. These observations indicate that not all the available reserves in the system are valuable from the system operator's point of view.

In many countries, electricity markets are still evolving with the challenges of uncertain energy resource and load [19]. The ramping products are proposed by California ISO [20] to accommodate the uncertainties. It should be emphasized that the bus-level delivery is not considered for ramping capacity in [20]. Before applying the robust optimization SCED in the real market, the corresponding pricing theory is imperative. This paper tries to propose some new ideas to clear this obstacle. The three major contributions of this paper are listed as follows.

1) A new concept, deliverable generation ramping reserve, is proposed within the AARO SCED framework. The generation ramping reserve is the additional ramping capability of the generator when part of it is "locked" in the SCED schedule.
2) The prices for the generation ramping reserves as well generation capacity reserves are derived within the robust co-optimization framework. With the help of the price information, the valuable reserves can be easily identified among the available reserves.
3) The market equilibrium is characterized by the proposed prices and dispatch instructions. The market participants can get the maximal profit by following ISO's the dispatch instruction and price signals.

The rest of this paper is organized as follows. The derivation of the reserve prices is presented in Section II based on AARO SCED, and then the market mechanism to credit flexibilities is discussed. Case studies for 3-Bus and IEEE 118-Bus systems are presented in Section III. Section IV concludes this paper.

## II. AARO SCED AND PRICES

In electricity markets, RTOs/ISOs normally operate two markets including DAM and RTM (or balancing market) [21]. The majority of the trades is cleared in DAM via SCUC [21], and SCED is normally performed periodically in RTM. This paper mainly focuses on the RTM.

In standard SCED problem, the unit generation output is subject to the following constraints which include unit capacity limits (1) and unit ramping up/down limits (2) (3).

$$\hat{I}_{i,t} P_i^{\min} \leq P_{i,t} \leq \hat{I}_{i,t} P_i^{\max}, \forall i, t \quad (1)$$
$$P_{i,t} - P_{i,(t-1)} \leq R_i^U (1 - \hat{y}_{i,t}) + P_i^{\min} \hat{y}_{i,t}, \forall i, t \quad (2)$$
$$-P_{i,t} + P_{i,(t-1)} \leq R_i^D (1 - \hat{z}_{i,t}) + P_i^{\min} \hat{z}_{i,t}, \forall i, t \quad (3)$$

Equations (2) and (3) show that a unit has to operate at its minimum capacity in two cases: right after it is turned on or right before it is turned off, which implies that the unit cannot provide reserve in those two cases. For notation brevity, we use matrix and vector to replace the formulations above. Then the SCED problem can be formulated as

$$\min_{\boldsymbol{P_i}} \quad \sum_i \boldsymbol{P_i}^\top \boldsymbol{Q_i} \boldsymbol{P_i} + \boldsymbol{B_i}^\top \boldsymbol{P_i} + f_i \quad (4)$$

$$\text{s.t.} \quad (\boldsymbol{\lambda}) \quad \sum_i \boldsymbol{P_i} = \boldsymbol{Dd}, \quad (5)$$

$$(\boldsymbol{\alpha_i} \geq \boldsymbol{0}) \quad \boldsymbol{AP_i} \leq \boldsymbol{R_i}, \forall i \quad (6)$$

$$(\boldsymbol{\eta} \geq \boldsymbol{0}) \quad \sum_i \boldsymbol{\Gamma_i} \boldsymbol{P_i} - \boldsymbol{\Gamma_d} \boldsymbol{d} \leq \boldsymbol{F}, \quad (7)$$

where (4) stands the operation cost. (5) denotes the load balance constraints. (6) is the compact form of (1)-(3). (7) represent the transmission constraints. $\boldsymbol{\Gamma}_i \in \mathbb{R}^{2N_L N_T \times N_T}$, and $\boldsymbol{\Gamma}_d \in \mathbb{R}^{2N_L N_T \times N_D N_T}$.

Due to the forecasting errors of renewable power output and load, ISOs/RTOs need to run SCED on a rolling basis in real time to balance the system. Recently, some approaches have been successfully applied in the SCUC/SCED problem to address the uncertainty issues caused by variations of load and renewables. Both stochastic SCUC/SCED and robust SCUC/SCED are studied intensively when considering the uncertainties.

To our best knowledge, this paper represents the first work on pricing the reserves in robust optimization framework. Hence, the following assumptions are made so that we can focus on the concept.

- Transmission loss is ignored in the SCED problem.
- The proposed approach is for ex ante dispatch and ex ante pricing. It is assumed that units are dispatched according to these instructions during the scheduling intervals.
- Units bid only energy price. The reserve bid is zero.
- Uncertainty comes from loads. Renewables are treated as negative loads. Other uncertainties such as contingencies are not discussed in this paper.
- Uncertainty set information is available to the ISO/RTO.
- Expectation and covariance of uncertainties can be obtained from historical data. However, PDF information is not specified for uncertainty as it is hard to obtain for uncertain sources such as RES.

### A. Affinely Adjustable Robust SCED

The basic idea of AARO optimization is originally from paper [15], in which a linear "feedback" in control theory is used to adjust dispatch with the realization of load. Authors in [13], [14], [22] applied it to solve SCED problem. In this section, the generation output is affinely adjusted according to the uncertainties,

$$\hat{\boldsymbol{P}}_i = \boldsymbol{P}_i + \boldsymbol{G}_i \boldsymbol{\epsilon}, \forall i, \tag{8}$$

where $\boldsymbol{G}_i \in \mathbb{R}^{N_T \times N_D N_T}$ is the affine adjustment matrix. $\boldsymbol{P}_i \in \mathbb{R}^{N_T}$ and $\hat{\boldsymbol{P}}_i \in \mathbb{R}^{N_T}$ are base dispatch and adjusted dispatch respectively, and $\boldsymbol{\epsilon} \in \mathbb{R}^{N_D N_T}$ is the uncertainty vector (i.e., deviation of loads from forecasted values). The new unit dispatch can be regulated based on the load deviation. It is noted that $\boldsymbol{\epsilon} \in \mathcal{U}$, and

$$\mathcal{U} := \{\boldsymbol{\epsilon} : \boldsymbol{S}\boldsymbol{\epsilon} \leq \boldsymbol{h}\}, \tag{9}$$

where $\boldsymbol{S} \in \mathbb{R}^{N_K \times N_D N_T}$ and $\boldsymbol{h} \in \mathbb{R}^{N_K}$. $\mathcal{U}$ is a polyhedron which includes more than just the lower and upper bounds for uncertainty. It is noted that entry in $\boldsymbol{h}$ is considered as uncertainty level which is positive. The AARO SCED, denoted as (ROP), can be formulated as

$$(\text{ROP}) : \min_{\boldsymbol{P}_i, \boldsymbol{G}_i} \sum_i \mathbb{E}\left[\hat{\boldsymbol{P}}_i^\top \boldsymbol{Q}_i \hat{\boldsymbol{P}}_i + \boldsymbol{B}_i^\top \hat{\boldsymbol{P}}_i + f_i\right] \tag{10}$$

$$\text{s.t.} \quad \sum_i \hat{\boldsymbol{P}}_i = \boldsymbol{D}(\boldsymbol{d} + \boldsymbol{\epsilon}), \forall \boldsymbol{\epsilon} \in \mathcal{U} \tag{11}$$

$$\boldsymbol{A}\hat{\boldsymbol{P}}_i \leq \boldsymbol{R}_i, \forall i, \forall \boldsymbol{\epsilon} \in \mathcal{U} \tag{12}$$

$$\sum_i \boldsymbol{\Gamma}_i \hat{\boldsymbol{P}}_i - \boldsymbol{\Gamma}_d (\boldsymbol{d} + \boldsymbol{\epsilon}) \leq \boldsymbol{F}, \forall \boldsymbol{\epsilon} \in \mathcal{U} \tag{13}$$

(8),

where $\boldsymbol{Q}_i$ is a semi-definite matrix with only diagonal entry. The objective function (10) is to minimize the expected cost. Without loss of generality, assume $\mathbb{E}(\boldsymbol{\epsilon}) = \boldsymbol{0}$. The objective function (10) can be rewritten as

$$\sum_i \left(\boldsymbol{P}_i^\top \boldsymbol{Q}_i \boldsymbol{P}_i + \boldsymbol{B}_i^\top \boldsymbol{P}_i + \Lambda(\boldsymbol{G}_i)\right),$$

where $\Lambda(\boldsymbol{G}_i) = \text{Tr}[\boldsymbol{G}_i^\top \boldsymbol{Q}_i \boldsymbol{G}_i \mathbb{E}(\boldsymbol{\epsilon}\boldsymbol{\epsilon}^\top)] + f_i$. It is assumed that covariance matrix $\mathbb{E}(\boldsymbol{\epsilon}\boldsymbol{\epsilon}^\top)$ is available, but PDF is unavailable. The problem (ROP) is converted to a computationally tractable problem (P) as follows, where the constraints including uncertain parameters are exactly reformulated.

$$(\text{P}) : \min_{\boldsymbol{P}_i, \boldsymbol{G}_i, \boldsymbol{\rho}_i, \boldsymbol{\zeta}} \sum_i \left(\boldsymbol{P}_i^\top \boldsymbol{Q}_i \boldsymbol{P}_i + \boldsymbol{B}_i^\top \boldsymbol{P}_i + \Lambda(\boldsymbol{G}_i)\right) \tag{14}$$

$$\text{s.t.} \quad (\boldsymbol{\lambda}) \quad \sum_i \boldsymbol{P}_i = \boldsymbol{D}\boldsymbol{d} \tag{15}$$

$$(\boldsymbol{\gamma}) \quad \sum_i \boldsymbol{G}_i = \boldsymbol{D} \tag{16}$$

$$(\boldsymbol{\alpha}_i \geq \boldsymbol{0}) \quad -\boldsymbol{R}_i + \boldsymbol{A}\boldsymbol{P}_i + \boldsymbol{\rho}_i^\top \boldsymbol{h} \leq \boldsymbol{0}, \forall i \tag{17}$$

$$(\boldsymbol{\beta}_i) \quad \boldsymbol{A}\boldsymbol{G}_i - \boldsymbol{\rho}_i^\top \boldsymbol{S} = \boldsymbol{0}, \forall i \tag{18}$$

$$(\boldsymbol{\eta} \geq \boldsymbol{0}) \quad -\boldsymbol{F} - \boldsymbol{\Gamma}_d \boldsymbol{d} + \sum_i \boldsymbol{\Gamma}_i \boldsymbol{P}_i + \boldsymbol{\zeta}^\top \boldsymbol{h} \leq \boldsymbol{0} \tag{19}$$

$$(\boldsymbol{\tau}) \quad \sum_i \boldsymbol{\Gamma}_i \boldsymbol{G}_i - \boldsymbol{\Gamma}_d - \boldsymbol{\zeta}^\top \boldsymbol{S} = \boldsymbol{0} \tag{20}$$

$$\boldsymbol{\rho}_i \geq \boldsymbol{0}, \boldsymbol{\zeta} \geq \boldsymbol{0},$$

where $\boldsymbol{\rho}_i^\top \in \mathbb{R}^{4N_T \times N_K}$ and $\boldsymbol{\zeta}^\top \in \mathbb{R}^{2N_L N_T \times N_K}$ are also variables. (15) and (16) are derived from (8) and (11). (17)-(20) are obtained from strong duality. Problem (P) is *convex* and can be solved efficiently by commercial solvers such as CPLEX and GUROBI. It can be observed that no PDF information is required to solve (P).

Different from the standard SCED (4)-(7), extra terms $\boldsymbol{\rho}_i^\top \boldsymbol{h}$ and $\boldsymbol{\zeta}^\top \boldsymbol{h}$ are added in inequality constraints (17) and (19), respectively, in problem (P). As $\boldsymbol{\rho}_i^\top \geq \boldsymbol{0}, \boldsymbol{\zeta}^\top \geq \boldsymbol{0}$ and $\boldsymbol{h} \geq \boldsymbol{0}$, $\boldsymbol{\rho}_i^\top \boldsymbol{h}$ and $\boldsymbol{\zeta}^\top \boldsymbol{h}$ are non-negative. It indicates that certain unit and transmission constraints in standard SCED are replaced with stronger constraints in the robust framework. Thus, the system actually keeps certain flexibilities for uncertainty accommodation. In the following sections, we will analyze these flexibilities, which are also called reserves.

### B. Generation Ramping Reserve and Capacity Reserve

An Uncertainty Mitigator (UM) refers to a flexible resource provider that participates in the management of uncertainties.

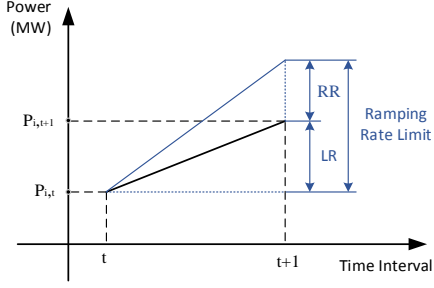

Fig. 1. Illustration of Upward Ramping Reserve (LR: Locked Ramping, RR: Ramping Reserve)

The flexible resources include generators with available ramping capabilities and adjustable loads. UMs have to keep certain reserves in order to accommodate the uncertainties. Compared to (6), constraint (17) may be binding even when the scheduled dispatch does not reach the capacity limits (1) or ramping limits (2)-(3). Based on the optimal solution to (P), the reserves can be calculated.

The generation ramping reserve is defined as the unused unit ramping capability or the value of the slack variable for (2-3), i.e.

$$Q_{i,t}^{ru} := R_i^U(1 - \hat{y}_{i,t}) + \frac{P_i^{\min}\hat{y}_{i,t} - (P_{i,t}^* - P_{i,(t-1)}^*)}{t - (t-1)}, \forall i, t, \quad (21)$$

$$Q_{i,t}^{rd} := R_i^D(1 - \hat{z}_{i,t}) + \frac{P_i^{\min}\hat{z}_{i,t} - (P_{i,(t-1)}^* - P_{i,t}^*)}{t - (t-1)}, \forall i, t, \quad (22)$$

where $Q_{i,t}^{ru}$ and $Q_{i,t}^{rd}$ are the upward and downward ramping reserves, respectively. Fig. 1 illustrate the concept of the upward ramping reserve. With scheduled generation outputs at $t$ and $t+1$, the locked ramping or the used unit ramping capability is $(P_{i,t+1} - P_{i,t})/1$. The generation ramping reserve is the available ramping rate less the locked ramping. The ramping reserve guarantees that the uncertainty mitigator still has additional ramping capability during the scheduled ramping process. In the existing AS market, the "locked" ramping is normally ignored for the spinning reserves [23], [24]. In the RTM, the time resolution is 5 minutes or 15 minutes, if we ignore the ramping process, there is a chance that system cannot provide enough ramping capability to accommodate the uncertainties. In contrast, the time resolution is one hour in DAM, and the units have enough time to redispatch even if the locked ramping is ignored [17], [23].

The generation capacity reserve is defined as unused unit generation capacity or the constraint slack for (1), i.e.

$$Q_{i,t}^{cu} := \hat{I}_{i,t}P_i^{\max} - P_{i,t}^*, \quad (23)$$
$$Q_{i,t}^{cd} := P_{i,t}^* - \hat{I}_{i,t}P_i^{\min}, \forall i, t \quad (24)$$

where $Q_{i,t}^{cd}$ and $Q_{i,t}^{cu}$ are the lower and upper generation capacity reserves, respectively. From the system's point of view, the total capacity reserves are fixed when the unit commitment and load level are determined.

### C. Marginal Prices

In most generation and reserve co-optimization approaches in AS market, the explicit reserve requirement constraint is modeled [17], [18], [25]. The shadow prices of this type of constraint is employed to derive the reserve price which reflects the coupled effects of the generation and reserve. Instead of setting the reserve manually and heuristically based on Monte Carlo simulations, the reserves in the AARO SCED are determined automatically in one shot. Although it has obvious advantages over the traditional reserve determination, it also poses new challenges on reserve price derivation. The existing pricing approaches cannot be used directly [17], [18] due to the lack of explicit reserve requirement constraints.

While the amount of reserves can be calculated according to (21)-(24), the question is how to set the prices for them. On the one hand, it is well known that not all the reserves are deliverable for uncertainty accommodation due to network constraints. On the other hand, the generation and reserve are coupled together in the RTM even if the reserve bid price is zero [26], and the market clearing price for reserve has certain relations with that for energy. In this paper, we call the reserve obtained by (21)-(24) *available reserve*. If a small increment or decrement of the reserve amount causes change of the expected operation cost in (P), then this type of reserve is called *valuable reserve*. To determine the exact value a reserve has, we derive the marginal prices for the reserve according to the Lagrangian function as follows.

The Lagrangian function for (P) is

$$\Phi(\boldsymbol{P_i}, \boldsymbol{G_i}, \boldsymbol{\rho_i}, \boldsymbol{\zeta}, \boldsymbol{\lambda}, \boldsymbol{\gamma}, \boldsymbol{\alpha_i}, \boldsymbol{\beta_i}, \boldsymbol{\eta}, \boldsymbol{\tau}) \quad (25)$$
$$= \sum_i \left(\boldsymbol{P_i}^\top \boldsymbol{Q_i} \boldsymbol{P_i} + \boldsymbol{B_i}^\top \boldsymbol{P_i} + \Lambda(\boldsymbol{G_i})\right)$$
$$+ \boldsymbol{\lambda}^\top (\boldsymbol{Dd} - \sum_i \boldsymbol{P_i}) + \text{Tr}[\boldsymbol{\gamma}^\top (\boldsymbol{D} - \sum_i \boldsymbol{G_i})]$$
$$+ \sum_i \left(\boldsymbol{\alpha_i}^\top (\boldsymbol{AP_i} + \boldsymbol{\rho_i}^\top \boldsymbol{h} - \boldsymbol{R_i}) + \text{Tr}[\boldsymbol{\beta_i}^\top (\boldsymbol{AG_i} - \boldsymbol{\rho_i}^\top \boldsymbol{S})]\right)$$
$$+ \boldsymbol{\eta}^\top \left(-\boldsymbol{F} - \boldsymbol{\Gamma_d d} + \sum_i \boldsymbol{\Gamma_i P_i} + \boldsymbol{\zeta}^\top \boldsymbol{h}\right)$$
$$+ \text{Tr}[\boldsymbol{\tau}^\top (\sum_i \boldsymbol{\Gamma_i G_i} - \boldsymbol{\Gamma_d} - \boldsymbol{\zeta}^\top \boldsymbol{S})]$$

The ramping reserve price for UM $i$ is defined as the marginal expected cost due to a unit decrement of the generation ramping rate of UM $i$. The capacity reserve price for UM $i$ is defined as the marginal expected cost due to a unit decrement of the generation capacity of UM $i$. They can be obtained from the Lagrangian function as

$$\boldsymbol{\pi_i^r} = -\frac{\partial \Phi}{\partial \boldsymbol{R_i}} = \boldsymbol{\alpha_i^*}. \quad (26)$$

Hence, the reserve provided by $i$ is valuable only when $\boldsymbol{\pi_i^r}$ is non-zero. Note that $\boldsymbol{\alpha_i^*}$ consists of Lagrangian multipliers for unit ramping limits and unit capacity limits.

The LMP for AARO SCED can be obtained based on its definition. It is the marginal expected cost due to a unit increment of the load. For generator $i$, it is formulated as

$$\boldsymbol{\pi_i^e} = \boldsymbol{\lambda^*} - \boldsymbol{\Gamma_i}^\top \boldsymbol{\eta^*}, \quad (27)$$

where $\pi_i^e \in \mathbb{R}^{N_T}$. $\pi_i^e$ also consists of energy component and congestion component. As shown in [13], the affine adjustment price can be obtained as

$$\pi_i^a = \gamma^* - \Gamma_i^\top \tau^* - A^\top \beta_i^*, \quad (28)$$

$\pi_i^a \in \mathbb{R}^{N_T \times N_D N_T}$ represents the marginal value of the adjustment coefficient in $G_i$. As stated in [13], $\text{Tr}[G_i^\top \pi_i^a]$ is the payment to $i$. It is noted that the unit of $\pi_i^a$ is \$ as $G_i$ is participant factor.

### D. Credit to Uncertainty Mitigators

Within the AARO framework, UMs help the system withstand the "load deviation" in the future. As the generation and the reserve are coupled together, the credit to reserve should reflect the coupling effect. Only the UMs who provide the valuable reserves are entitled to credits.

The price of the flexible resources are $\pi_i^r$ for UM $i$. The total credit allocated to $i$ is

$$\Theta_i = (R_i - AP_i^*)^\top \pi_i^r, \quad (29)$$

which is the product of the reserve price and the reserve amount. In fact, the reserve price $\pi_i^r$ reflects how much "value" the reserve has. $R_i - AP_i^*$ reflects the available reserve, which is the reserve quantity at each time interval. Only when the reserve is a valuable reserve (i.e. $\pi_i^r \neq 0$), UM $i$ get the reserve credit. Otherwise, the credit entitled to UM $i$ is zero even if the available reserve it provides is non-zero. There are similar phenomenons in the traditional zonal-based reserve market. For example, the reserve price \$0/MWh at a specific zone occurs when the cleared system reserve is higher than the required amount in Case I of Section VI [17].

### E. Market Equilibrium

In the partial market equilibrium model, it is assumed that market participants are price takers [27]. This assumption is popular in the electricity market [25], [28]. The expected profit maximization problem for source $i$ can be formulated as

$$(\text{PMP}_i): \max_{P_i, G_i} \mathbb{E} \left\{ \begin{array}{l} P_i^\top \pi_i^e + \text{Tr}[G_i^\top \pi_i^a] - B_i^\top (P_i + G_i \epsilon) \\ -(P_i + G_i \epsilon)^\top Q_i (P_i + G_i \epsilon) - f_i + \Theta_i \end{array} \right\}$$
$$\text{s.t.} \quad (8)(12)$$

where the decision variables for source $i$ are $P_i$ and $G_i$, given the price signal $\pi_i^e, \pi_i^a$, and $\pi_i^r$. Credit $\text{Tr}[G_i^\top \pi_i^a]$ is for the affine adjustment matrix as stated in [13]. It should be noted that the constraint for mitigator $i$ is (12), which is why incentive $\Theta_i$ is necessary for $i$ to follow the dispatch instruction. The objective function can be converted to

$$P_i^\top \pi_i^e + \text{Tr}[G_i^\top \pi_i^a] + (R_i - AP_i)^\top \pi_i^r$$
$$-P_i^\top Q P_i - B_i^\top P_i - f_i - \text{Tr}\left[G_i^\top Q_i G_i \mathbb{E}(\epsilon \epsilon^\top)\right] \quad (30)$$

It can be observed that (30) is a portion of the Lagrangian function (25). The optimal solution $(P_i, G_i)$ to $(\text{PMP}_i)$ is a function of $(\pi_i^e, \pi_i^a, \pi_i^r)$. Since the problem (P) is convex and Slater's condition is satisfied, the strong duality holds. Therefore, the saddle point $(P_i^*, G_i^*)$, which is the optimal solution to (P), is also the optimal solution to $(\text{PMP}_i)$. Consequently, source $i$ can obtain the maximum profit by following the ISO's dispatch instruction $(P_i^*, G_i^*)$. If the objective function (14) is strictly convex, then $(P_i^*, G_i^*)$ is the unique choice for $i$ to get the maximal profit. Price signal $\pi_i^r$ (associated with credit $\Theta_i$) and $\pi_i^e$ provide the incentives for $i$ to dispatch power output to $P_i^*$, which supplies the load and maintains the ramping and capacity reserves. In addition, the price signal $\pi_i^a$ provides the incentive for $i$ to follow the adjustment instruction $G_i$.

Given that the optimality condition of problem (P) is satisfied, dispatch signal $(P_i^*, G_i^*)$ and price signal $(\pi_i^e, \pi_i^a, \alpha_i^*)$ constitute a competitive partial equilibrium [27]. The assumption made for the equilibrium model is that the market participant uses affine policy to adjust the generation, so the dispatch obtained from AARO SCED is near-optimal solution. Consequently, the equilibrium conditions are also near-optimal for market participants.

An alternative way to get the partial market equilibrium is to lump the reserve price into the LMP. However, it has serious incentive issues. The integrated LMP can be written as

$$\hat{\pi}_i^e = \lambda^* - \Gamma_i^\top \eta^* - A^\top \alpha_i. \quad (31)$$

The new credit UM $i$ receives on the energy and reserve is

$$(\hat{\pi}_i^e)^\top P_i = (\pi_i^e)^\top P_i - (\alpha_i)^\top AP_i. \quad (32)$$

It can be observed that the reserve credit is negative when the shadow price of the upper bound constraint is nonzero. It means that the more flexible resources UM provides, the fewer profit it gets. The UM is provided negative incentives for the flexible resources. In contrast, the reserve credit defined in (29) is always nonnegative.

## III. CASE STUDY

A 3-Bus system and the modified IEEE 118-bus system are studied in this section to illustrate the concepts of available/valuable ramping/capacity reserves and the associated prices, as well as their impacts on market participants. The simulations are carried out by CPLEX 12.5 on PC with Intel i7-3770 3.40GHZ 8GB RAM.

### A. 3-Bus System

The 3-Bus system is consisted of two units, one wind farm, two loads, and three lines. Please refer http://motor.ece.iit.edu/data/rscuc/3_bus_data.pdf for the detailed data. For simplicity, three time intervals are studied and a single-segment incremental cost (IC) is employed to represent the fuel cost. The time resolution is 15 minutes. It is assumed that both units in the system are committed, which is determined in the DAM. Load in the current interval (i.e., Time 1) is assumed to be accurate. Loads from Time 2 and Time 3 are forecasted based on current available information, and forecasting errors may exist. Assume that the expectation of uncertainties at Time 2 and Time 3 are 0, and their probability distribution is unknown.

The total expected cost calculated based on the AARO SCED is \$1380.9. It is higher than the standard SCED cost \$1333.5, which cannot be immuned against uncertainties. It indicates that the expensive G2 supplies more loads within AARO SCED than it does in standard SCED. The reason

TABLE I
BASE ED AND LMP FOR 3-BUS SYSTEM IN CASE 2

| Time | G1 (MW) | G2 (MW) | LMP ($/MWh) | | |
|---|---|---|---|---|---|
| | | | Bus1 | Bus2 | Bus3 |
| 1 | 126.5 | 13.5 | 10 | 10 | 10 |
| 2 | 131.65 | 23.35 | 10 | 32.5 | 23.5 |
| 3 | 143.6 | 23.4 | 10 | 32.5 | 23.5 |

TABLE II
UPWARD RESERVES AND PRICES FOR 3-BUS SYSTEM IN CASE 2

| Time | G1 | | | | G2 | | | |
|---|---|---|---|---|---|---|---|---|
| | RR | RRP | CR | CRP | RR | RRP | CR | CRP |
| 1 | 18.5 | 0 | 53.5 | 0 | 6.5 | 0 | 66.5 | 0 |
| 2 | 19.85 | 0 | 48.35 | 0 | 0.15 | 15 | 56.65 | 0 |
| 3 | 13.05 | 0 | 36.4 | 0 | 9.95 | 7.5 | 56.6 | 0 |

RR: Ramping Reserve, MW/15min; RRP: Ramping Reserve Price, $/(MWh/15min)
CR: Capacity Reserve, MW; CRP: Capacity Reserve Price, $/MWh

TABLE III
CREDITS ENTITLED TO UNITS BASED ON RESERVE PRICES IN CASE 2 ($)

| Time | G1 | G2 |
|---|---|---|
| 1 | 0 | 0 |
| 2 | 0 | 15*0.15*0.25=0.5625 |
| 3 | 0 | 7.5*9.95*0.25=18.6563 |

of dispatching the extra generation from G2 is that operator needs additional deliverable reserves to accommodate the uncertainties while minimizing the total cost. The non-zero entries in $G_1 \in \mathbb{R}^{3 \times 6}$ are $g_{11}^1 = g_{12}^2 = 1, g_{23}^1 = 1, g_{24}^1 = 0.98, g_{36}^1 = 0.4$, and the non-zero entries in $G_2 \in \mathbb{R}^{3 \times 6}$ are $g_{24}^2 = 0.02, g_{35}^2 = 1, g_{36}^2 = 0.6$. If the loads on bus 1 and bus 3 at interval 2 are increased to 97MW and 67.5MW, respectively, then the units can be re-dispatched based on the adjustment matrix to $131.65 + 1 \times 7 + 0.98 \times 2.5 = 141.025$MW and $23.35 + 0 \times 7 + 0.02 \times 2.5 = 23.475$MW.

The question is whether the units have incentives to maintain the reserves. First, we consider the scenario without reserve credits. The LMP at bus 2 at Time 2 of $32.5/MWh is larger than G2's marginal cost of $10/MWh. Therefore, G2 is inclined to supply more loads at Time 2. The increase of the output from G2 will shrink the reserve it can provide. The fact that G2 has negative profit at Time 1 is the uplift issue [29], which is beyond the scope of this paper. In this scenario, the market participants would ignore the uncertainties and game the market.

Now, we consider the second scenario with reserve credits. The upward reserves provided by UMs are shown in Table II. It is observed that the capacity reserve is not scarce resource, i.e. the online capacity is adequate. They are just available reserves. In contrast, the prices of upward ramping reserve at Time 2 and Time 3 are non-zeros. According to the definition in this paper, they are valuable reserves. In this scenario, UMs are entitled certain credits based on the contribution. As shown in Table III, uncertainty mitigator G2 is entitled $0.56 at Time 2 and $18.66 at Time 3, where $15/(MWh/15min) and $7.5/(MWh/15min) are the entries in $\pi_2^r$ for ramping constraints at Time 2 and Time 3, respectively; and 0.15

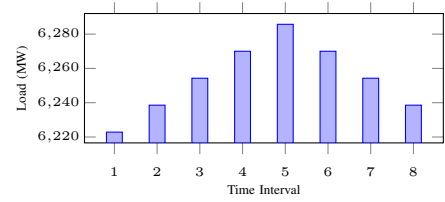

Fig. 2. System-wide Net Load

MW/15min and 9.95 MW/15min are the ramping reserves of G2 at Times 2 and 3, respectively. The credit is calculated based on (29). It is noted that the credits are time-coupled, as the ramping reserve is calculated according to the difference of the power outputs between two intervals.

Next we give an example to show whether G2 is inclined to deviate from the base dispatch instruction with the credits from LMPs and $\pi_i^r$. If G2 follows the instruction $P_2^*$ as shown in Table I, then the total credit associated with LMPs and $\pi_i^r$ are $[13.5, 23.35, 23.4] \cdot [10, 32.5, 32.5]^\top \cdot 0.25 + 0.5625 + 18.6563 = \$432.8125$. Consider a possible base dispatch $[12.5, 22.35, 23.4]$, where G2 generates 1MW less at Time 1 as its profit is negative at Time 1, and it also generates 1MW less at Time 2 to maintain the same ramping reserve. Then the new credit is $[12.5, 22.35, 23.4] \cdot [10, 32.5, 32.5]^\top \cdot 0.25 + 0.5625 + 7.5 \cdot 8.95 \cdot 0.25 = \$420.3125$, where ramping reserve at Time 3 is decreased to 8.95MW/15min. The credit is $12.5 lower than that by following $P_2^*$, and the saved fuel cost by supplying smaller load is $2 \cdot 25 \cdot 0.25 = \$12.5$. Therefore, G2 does not get more profit by using this dispatch. This example illustrates how G2 gets credit associated with LMPs and $\pi_i^r$, and why G2 is not inclined to deviate the dispatch instruction $P_2^*$. The more rigorous mathematical analysis on market equilibrium is shown in Section II-E.

### B. Modified IEEE 118-Bus System

There are 54 traditional units and 186 branches in the modified IEEE 118-Bus system. The scheduling period is 2 hours, and the time interval is 15 minutes. The loads are depicted in Fig. 2. The UCs are determined in advance by the solution to robust SCUC problem with 5% reserves. Five wind farms are introduced in the system, and they are located at buses 11, 49, 60, 78, and 90, respectively. We denote the set of buses with uncertainty as $\mathcal{M}$. It is assumed that the forecasted power output (i.e. nominal output) and installed capacity for each wind farm are 100 MW and 200 MW, respectively. The uncertainties in this case are from the RES only. The uncertainty $\epsilon_{m,t}$ satisfies

$$\begin{cases} |\epsilon_{m,t}| \leq 100 \cdot r_1 \big(1 + 0.01 \cdot (t-1)\big), m \in \mathcal{M}, \forall t \\ \Big|\sum_m \epsilon_{m,t}\Big| \leq 500 \cdot r_1 \cdot r_2 \big(1 + 0.01 \cdot (t-1)\big), m \in \mathcal{M}, \forall t \end{cases}$$

where $r_1$ reflects the forecast error confidence interval for a single wind farm [4], [14], and $r_2$ reflects the forecast error confidence interval for the aggregated wind output. When $r_2 < 1$, it indicates that the aggregated forecast error confidence interval is smaller than the sum of five single intervals. In the experiment, the forecast error increases with the time intervals.

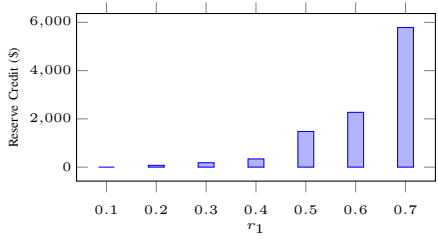

Fig. 3. Total Reserve Credit W.R.T. Uncertainty Levels and Fixed Nominal Wind Power ($r_2 = 1$)

TABLE IV
RESERVES W.R.T INCREASING UNCERTAINTY LEVELS AND FIXED NOMINAL WIND POWER ($r_2 = 1, t = 2$)

| $r_1$ | Upward Reserve | | | | Downward Reserve | | | |
|---|---|---|---|---|---|---|---|---|
| | AvaRamp | ValRamp | AvaCap | ValCap | AvaRamp | ValRamp | AvaCap | ValCap |
| 0.1 | 886.79 | 0 | 661.43 | 0 | 938.46 | 0 | 4203.57 | 0 |
| 0.2 | 886.79 | 180.79 | 661.43 | 10.61 | 938.46 | 10.37 | 4203.57 | 0 |
| 0.3 | 886.79 | 180.29 | 661.43 | 61.36 | 938.46 | 0 | 4203.57 | 0 |
| 0.4 | 886.79 | 179.79 | 661.43 | 112.11 | 938.46 | 0 | 4203.57 | 0 |
| 0.5 | 886.79 | 479.29 | 661.43 | 162.86 | 938.46 | 361.61 | 4203.57 | 0 |
| 0.6 | 886.79 | 513.33 | 661.43 | 215.49 | 938.46 | 337.5 | 4203.57 | 4.52 |
| 0.7 | 886.79 | 684.04 | 661.43 | 316.72 | 938.46 | 680.33 | 4203.57 | 12.14 |

AvaRamp: Available Ramping Reserve (MW/15min); ValRamp: Valuable Ramping Reserve (MW/15min);
AvaCap: Available Capacity Reserve (MW); ValCap: Valuable Capacity Reserve (MW);

The detailed data including unit parameters, uncertainty correlation matrix, line reactance and ratings, and net load profiles can be found at http://motor.ece.iit.edu/Data/118_UMP.xls.

We consider the interval bounds for the uncertainties, and perform the sensitivity analysis with respect to $r_1$. Fig. 3 shows the reserve credits (RC) the UMs receive with the change of the $r_1$. The reserve credit (RC) is the sum of the ramping reserve credit and capacity reserve credit UMs are entitled to. When $r_1$ is high, the UMs are also entitled to high credits. As shown in (29), the reserve credits are the sum of the products of the amount of valuable reserve and the price of the valuable reserve. They are analyzed as follows.

Table IV presents the upward/downward available reserves and valuable reserves at Time 2 with increasing forecast errors, fixed normal wind power output (100MW each), and fixed $r_2$. The $r_1 \in [0.1, 0.7]$, so the error in percentage of the installed capacity (200MW each) is from 5% to 35%. It is observed that the available reserves remain the same while valuable reserve change dramatically with the forecast errors. As shown in Table IV, upward available ramping reserve remains 886.70MW/15min, and the capacity reserve remains 661.43MW. The main reason is that the unit commitment and load demand is fixed at Time 2 in the system. In contrast, the upward valuable reserve is 0 when $r_1 = 0.1$. It indicates that the opportunity cost of keeping the ramping reserve for UM is zero as it can recover the profit from the energy credit. When $r_1$ is 0.2, 0.3 and 0.4, the valuable ramping reserve is around 180MW/15min, and the valuable capacity reserve is around 11MW, 61MW, and 112MW respectively. UMs are entitled to credits by keeping the reserves, which is also shown in Fig. 3. It suggests that opportunity cost of keeping the reserve is non-zero, i.e., UMs can get more profits by deviating from the dispatch instruction if they are not entitled to reserve credits. When $r_1$ is further increased to 0.5, the amount of valuable ramping reserve jumps to 479MW/15min. It means that more available reserves become valuable when the uncertainty level is high. The valuable capacity reserve also increases to 163 MW in this case. A similar tendency can also be observed for the downward reserves shown in Table IV. It should be emphasized that the amount of valuable reserve doesn't change monotonically with the uncertainty level. Instead, what we revealed in this paper is a trend.

The upward available ramping reserve and the price of unit 42 are depicted for different time intervals in Fig. 4 with $r_1 = 0.7, r_2 = 1$. As shown in Fig. 4, at Time 7, the unit only provides available reserve but not valuable reserve as its price is zero. It can be observed that the ramping reserve price reaches its highest point at Time 5, which is also the peak load interval. In contrast, the ramping reserve price is low at Time 6 although the load demand at Time 6 is still relatively high compared to those at other intervals. It is observed the load climbs from Time 4 to Time 5, but falls from Time 5 to Time 6 as shown in Fig. 2. It indicates that the ramping reserve price is not only related with the load demand but also with the load change. In this case, the upward ramping reserve is scarce resource at Time 5, and the opportunity cost of keeping them is also high. In contrast, the upward ramping reserve is relatively cheap when the load demand is falling at Time 6, 7, and 8.

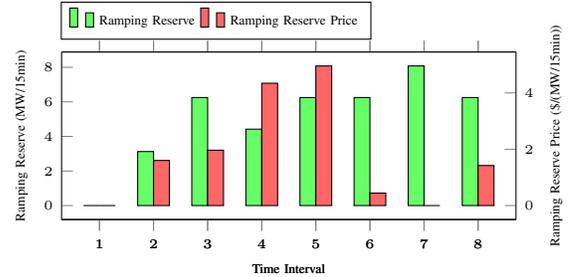

Fig. 4. Upward Ramping Reserve and Price for Unit 42 ($r_1 = 0.7, r_2 = 1$)

Fig. 5 depicts the available capacity reserve and the price of capacity reserve for unit 24 with $r_1 = 0.7, r_2 = 1$. Although the reserve amount unit 24 keeps is the same at each time interval, the price is different. It is observed that the capacity reserve price in this case has a similar trend to the system load level shown in Fig. 2. For example, the capacity reserve is the most expensive at Time 5 when the peak load occurs. The reason is that the system-wide upward capacity reserve is the online installed capacity, which is fixed, less the load level. When the load level is high, the upward capacity reserve is small which becomes a scarce resource in the system.

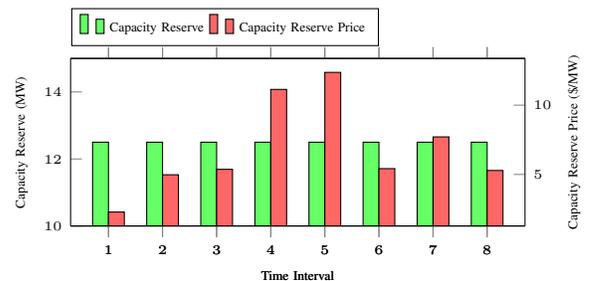

Fig. 5. Upward Capacity Reserve and Price for Unit 24 ($r_1 = 0.7, r_2 = 1$)

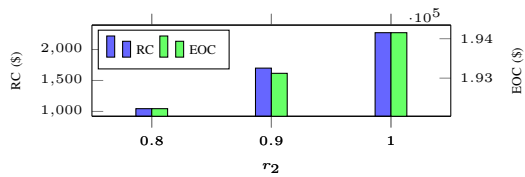

Fig. 6. Reserve Credit (RC) and Expected Operation Cost (EOC) with Different $r_2$ ($r_1 = 0.6$)

Reserve credits with respect to different $r_2$ are presented in Fig. 6. It shows that the decrease of $r_2$ (from 1 to 0.9 to 0.8) also leads to lower payments related to reserve. For example, when $r_2$ decreases from 1 to 0.8, the total RC decreases from around \$2,300 to \$1,000. The expected operation cost also decreases from \$194,000 to \$192,000. It indicates that the shrinking uncertainty set actually increases the feasible set for the robust dispatches.

The numerical results in this part indicate that the reserve payment proposed in this paper helps maximize the social warfare. When the uncertainty level is high, the payment related to reserve is also high, It may attract long-term investment of flexible resources. More flexible resources also mean the system has more capabilities to handle the uncertainties, and the system can accommodate higher RES penetrations.

## IV. Conclusions

This paper proposes a new concept, ramping reserve, within the AARO SCED framework. AARO SCED is an effective tool in RTM to address the uncertainty issue although its solution may only be near-optimal. The flexibilities in this paper include the generation ramping reserve and the generation capacity reserve. The prices for the ramping reserve and capacity reserve are also derived based on the Lagrangian function. They are the opportunity costs of the uncertainty mitigators to keep the reserves or flexibilities. With the help of these prices, the reserves are classified into two categories, which are available reserves and valuable reserves. The case studies explain the concept of these reserves and the impacts on the behaviors of market participants.

Many researches on this topic are open in the future. With increasing RES penetration in the power system, the flexibilities play a crucial role in uncertainty accommodation. The prices derived in this paper provide an option on how to provide the reserve signals within the robust optimization framework. Those reserve credits to UMs may attract the investment of flexible resources in the long term. In return, new flexibility investment allows the system to accommodate higher RES level.

It should be pointed out that the reserve prices are unit-specified. It is admitted that it is not perfect. It is ideal that all the resources at the same bus have the same price. An extension of the proposed reserve prices is to set the maximum ramping and capacity reserve prices for the units located at a node as the nodal price. In this way, the nodal reserve prices are determined by the most expensive opportunity cost for the reserve in the node. However, due to the affine policy and the non-zero optimality of the AARO SCED, it may deteriorate the partial market equilibrium to a certain extent.